\documentclass{article}
\usepackage{amsfonts}
\usepackage{amsmath}

\input{tcilatex}
\begin{document}

\begin{center}
\textbf{COMPARISON CRITERIA FOR DISCRETE FRACTIONAL STURM-LIOUVILLE EQUATIONS%
}
\end{center}

\bigskip

\begin{center}
Ramazan OZARSLAN$^{a\ast }$ and Erdal BAS$^{a}$

\bigskip

$^{a}$\textit{Firat University, Science Faculty, Department of Mathematics,
23119 Elazig/Turkey}

\bigskip

\textit{e-mail: }$^{\ast }$ozarslanramazan@firat.edu.tr\textit{,
erdalmat@yahoo.com}

\bigskip
\end{center}

\noindent \textbf{Abstract}. In this study, we give the Sturm comparison
theorems for discrete fractional Sturm-Liouville (DFSL) equations within
Riemann-Liouville and Gr\"{u}nwald-Letnikov sense. The emergence of
Sturm-Liouville equations began as one dimensional Schr\"{o}dinger equation
in quantum mechanics and one of the most important results is Sturm
comparison theorems \cite{levi}. These theorems give information about the
properties of zeros of two equations having different potentials.

\bigskip


\bigskip


\bigskip

\begin{center}
\textbf{1. Introduction}
\end{center}

Discrete fractional calculus has made great progress in recent years (see 
\cite{thabet,thabet2,anas2,ferhan,ferhan6,jin,okkes} ). It is discrete
fractional analogue of ODEs, for this reason some problems in ODEs theory
have been started to be applied similarly to discrete fractional problems.
Generally, this adaptation firstly was made to difference calculus, secondly
to fractional calculus and finally to discrete fractional calculus.
Fractional sums and differences were gained firstly in Diaz-Osler \cite{diaz}%
, Miller-Ross \cite{miller} and Gray and Zhang \cite{gray} and they obtained
various types of discrete fractional integrals and derivatives. Later,
several authors started to study on fractional difference calculus
Goodrich-Peterson \cite{goodrich}, Baleanu et.al. \cite{dumitru},
Mohan-Deekshitulu \cite{mohan}. However, discrete fractional calculus is a
rather novel area. The main studies has been done by Atici et.al. \cite%
{ferhan}, \cite{ferhan5}, Anastassiou \cite{anas2}, Abdeljawad et.al. \cite%
{thabet,thabet2,thabet4,thabet6} and Cheng et.al. \cite{jin} and so forth.

Discrete fractional calculus has still a lot of open problems when compared
to ODEs. One of the most important is Sturm-Liouville problems. Discrete
fractional Sturm-Liouville (DFSL) equation was firstly considered in \cite%
{almeida}, and next considered in our work \cite{ramo} and we showed the
fundemantal spectral properties, like self-adjointness of the operator,
orthogonality of the eigenfunctions, reality of the eigenvalues in \cite%
{ramo}. We considered DFSL equations in two different ways;

\noindent $i)$ (nabla right and left) Riemann-Liouville (\textbf{R-L)}
fractional operator,%
\begin{equation}
L_{1}x\left( t\right) =\nabla _{a}^{\mu }\left( p\left( t\right) _{b}\nabla
^{\mu }x\left( t\right) \right) +q\left( t\right) x\left( t\right) =\lambda
r\left( t\right) x\left( t\right) ,  \tag{1}  \label{a3}
\end{equation}

\noindent $ii)$ (delta left and right) Gr\"{u}nwald-Letnikov (\textbf{G-L)}
fractional operator,%
\begin{equation}
L_{2}x\left( t\right) =\Delta _{-}^{\mu }\left( p\left( t\right) \Delta
_{+}^{\mu }x\left( t\right) \right) +q\left( t\right) x\left( t\right)
=\lambda r\left( t\right) x\left( t\right) .  \tag{2}  \label{a4}
\end{equation}

The emergence of Sturm-Liouville equations began as one dimensional Schr\"{o}%
dinger equation in quantum mechanics and one of the most important results
is Sturm comparison theorems \cite{levi}. These theorems give information
about the properties of zeros of two equations having different potentials.
The first theorem asserts that 
\begin{eqnarray}
u^{\prime \prime }+g\left( x\right) u &=&0,  \TCItag{3}  \label{a1} \\
v^{\prime \prime }+h\left( x\right) v &=&0,  \TCItag{4}  \label{a2}
\end{eqnarray}%
$x\in \left[ a,b\right] ,$ if $g\left( x\right) <h\left( x\right) $ there is
at least one zero of each solution of $\left( \ref{a2}\right) $ between any
two zeros of any nontrivial solution of $\left( \ref{a1}\right) $.

The second theorem asserts that$,$ let $u\left( x\right) $ be the solution
of $\left( \ref{a1}\right) $ with the initial conditions 
\begin{equation*}
u\left( a\right) =\sin \alpha ,\text{ }u^{\prime }\left( a\right) =-\cos
\alpha ,
\end{equation*}%
and let $v(x)$ be the solution of $\left( \ref{a2}\right) $ with the same
initial conditions. Also, let $g(x)<h(x)$ for $x\in \left[ a,b\right] $.
Then if $u\left( x\right) $ has $m$ zeros in $a<x<b$, $v(x)$ has not fewer
than $m$ zeros in the same interval and the $k$th zero of $v(x)$ is less
than the $k$th zero of $u(x)$.

Various discrete versions of this theorem was considered in detailed in \cite%
{atkinson} and also was studied for infinite case in \cite{bohner} as%
\begin{eqnarray*}
\Delta ^{2}u_{k}+p_{k}u_{k+1} &=&0, \\
\Delta ^{2}v_{k}+p_{k}v_{k+1} &=&0.
\end{eqnarray*}%
\qquad Besides, fractional Sturm-Liouville differential operators have been
studied by \cite{bas,tuba,bas2,klimek}. More recently the fractional version
of this theorem has been considered in \cite{bas3} as, $0<\alpha <1$%
\begin{eqnarray*}
D_{b-}^{\alpha }p\left( r\right) ^{C}D_{a+}^{\alpha }u\left( r\right)
+g\left( r\right) u\left( r\right) &=&0, \\
D_{b-}^{\alpha }p\left( r\right) ^{C}D_{a+}^{\alpha }u\left( r\right)
+h\left( r\right) u\left( r\right) &=&0.
\end{eqnarray*}%
Additionally to these works, a different version of the comparison theorem
was considered in \cite{goodrich}$.$

In this work, discrete fractional versions of these theorems are discussed
by the equations $\left( \ref{a3}\right) $ and $\left( \ref{a4}\right) .$

\bigskip

\begin{center}
\textbf{2. Preliminaries}

\bigskip
\end{center}

\noindent \textbf{Definition 1. }\cite{bohner2} Falling factorial is defined
by, $\alpha \in 
\mathbb{R}
,$ 
\begin{equation*}
t^{\underline{\alpha }}=\frac{\Gamma \left( t+1\right) }{\Gamma \left(
t-\alpha +1\right) },
\end{equation*}%
where $\Gamma $ is the gamma function.

\bigskip

\noindent \textbf{Definition 2. }\cite{bohner2} Rising factorial is defined
by, $\alpha \in 
\mathbb{R}
,$ 
\begin{equation*}
t^{\bar{\alpha}}=\frac{\Gamma \left( t+\alpha \right) }{\Gamma \left(
t\right) }.
\end{equation*}%
Gamma function is well-defined in the above two definitions.

\bigskip

\noindent \textbf{Definition 3. }\cite{ferhan, thabet, miller} Fractional
sum operators are defined by,

\noindent $\left( i\right) $ The nabla left fractional sum of order $\mu >0$
is defined 
\begin{equation}
\nabla _{a}^{-\mu }x\left( t\right) =\frac{1}{\Gamma \left( \mu \right) }%
\sum_{s=a+1}^{t}\left( t-\rho \left( s\right) \right) ^{\overline{\mu -1}%
}x\left( s\right) ,\text{ }t\in 
\mathbb{N}
_{a+1},  \tag{5}  \label{6}
\end{equation}%
\noindent $\left( ii\right) $ The nabla right fractional sum of order $\mu
>0 $ is defined 
\begin{equation}
_{b}\nabla ^{-\mu }x\left( t\right) =\frac{1}{\Gamma \left( \mu \right) }%
\sum_{s=t}^{b-1}\left( s-\rho \left( t\right) \right) ^{\overline{\mu -1}%
}x\left( s\right) ,\text{ }t\in \text{ }_{b-1}%
\mathbb{N}
,  \tag{6}  \label{7}
\end{equation}%
where $\rho \left( t\right) =t-1$ is called backward jump operators, $%
\mathbb{N}
_{a}=\left \{ a,a+1,...\right \} ,$ $_{b}%
\mathbb{N}
=\left \{ b,b-1,...\right \} $.

\bigskip

\noindent \textbf{Definition 4. }\cite{thabet2, thabet4} Fractional
difference operators in \textbf{Riemann--Liouville (R-L) }sense are defined
by,

\noindent $\left( i\right) $ The nabla left fractional difference of order $%
\mu >0$ is defined 
\begin{equation}
\nabla _{a}^{\mu }x\left( t\right) =\nabla _{a}^{n}\nabla _{a}^{-\left(
n-\mu \right) }x\left( t\right) =\frac{\nabla ^{n}}{\Gamma \left( n-\mu
\right) }\sum_{s=a+1}^{t}\left( t-\rho \left( s\right) \right) ^{\overline{%
n-\mu -1}}x\left( s\right) ,\text{ }t\in 
\mathbb{N}
_{a+1},  \tag{7}  \label{9}
\end{equation}%
$\noindent \left( ii\right) $ The nabla right fractional difference of order 
$\mu >0$ is defined 
\begin{equation}
_{b}\nabla ^{\mu }x\left( t\right) =\left( -1\right) ^{n}\nabla
_{b}^{n}\nabla _{a}^{-\left( n-\mu \right) }x\left( t\right) =\frac{\left(
-1\right) ^{n}\Delta ^{n}}{\Gamma \left( n-\mu \right) }\sum_{s=a+1}^{t}%
\left( s-\rho \left( t\right) \right) ^{\overline{n-\mu -1}}x\left( s\right)
,\text{ }t\in \text{ }_{b-1}%
\mathbb{N}
.  \tag{8}  \label{10}
\end{equation}

\bigskip

\noindent \textbf{Definition 5. }\cite{bourdin,jin,diaz} Fractional
difference operators in \textbf{Gr\"{u}nwald--Letnikov (G-L) }sense are
defined by,

\noindent $\left( i\right) $ The delta left fractional difference of order $%
\mu ,$ $0<\mu \leq 1,$ is defined%
\begin{equation}
\Delta _{-}^{\mu }x\left( t\right) =\frac{1}{h^{\mu }}\sum_{s=0}^{t}\left(
-1\right) ^{s}\frac{\mu \left( \mu -1\right) ...\left( \mu -s+1\right) }{s!}%
x\left( t-s\right) ,\text{ }t=1,...,N.  \tag{9}  \label{11}
\end{equation}%
$\noindent \left( ii\right) $ The delta right fractional difference of order 
$\mu ,$ $0<\mu \leq 1,$ is defined%
\begin{equation}
\Delta _{+}^{\mu }x\left( t\right) =\frac{1}{h^{\mu }}\sum_{s=0}^{N-t}\left(
-1\right) ^{s}\frac{\mu \left( \mu -1\right) ...\left( \mu -s+1\right) }{s!}%
x\left( t+s\right) ,\text{ }t=0,..,N-1.  \tag{10}  \label{12}
\end{equation}

\bigskip

\noindent \textbf{Definition 6. }\cite{thabet4} Integration by parts formula
for\textbf{\ R-L} nabla fractional difference operator is defined by, $u$ is
defined on $_{b}%
\mathbb{N}
$ and $v$ is defined on $%
\mathbb{N}
_{a}$, 
\begin{equation}
\sum_{s=a+1}^{b-1}u\left( s\right) \nabla _{a}^{\mu }v\left( s\right)
=\sum_{s=a+1}^{b-1}v\left( s\right) _{b}\nabla ^{\mu }u\left( s\right) . 
\tag{11}  \label{16}
\end{equation}

\bigskip

\noindent \textbf{Definition 7. }\cite{almeida, bourdin} Integration by
parts formula for \textbf{G-L }delta fractional difference operator is
defined by, $u$, $v$ is defined on $\left \{ 0,1,...,n\right \} $,%
\begin{equation}
\sum_{s=0}^{n}u\left( s\right) \Delta _{-}^{\mu }v\left( s\right)
=\sum_{s=0}^{n}v\left( s\right) \Delta _{+}^{\mu }u\left( s\right) . 
\tag{12}  \label{17}
\end{equation}

\newpage

\begin{center}
\textbf{3. Main Results}

\bigskip
\end{center}

Let's consider the DFSL equation firstly in (\textbf{R-L)} sense; $0<\mu <1$%
\begin{eqnarray}
\nabla _{a}^{\mu }\left( p\left( t\right) _{b}\nabla ^{\mu }u\left( t\right)
\right) +k\left( t\right) u\left( t\right) &=&0,  \TCItag{13}  \label{c1} \\
\nabla _{a}^{\mu }\left( p\left( t\right) _{b}\nabla ^{\mu }v\left( t\right)
\right) +m\left( t\right) v\left( t\right) &=&0.  \TCItag{14}  \label{c2}
\end{eqnarray}

\bigskip

\noindent \textbf{Theorem 3.1. }(Sturm's 1st comparison theorem) Let $%
k\left( t\right) <m\left( t\right) $ over the entire interval $\left[ a,b%
\right] ,$ then between every two zeros of any nontrivial solution of the
equation $\left( \ref{c1}\right) $ there is at least one zero of every
solution of the equation $\left( \ref{c2}\right) $.

\bigskip

\textbf{Proof. }If we multiply the equation $\left( \ref{c1}\right) $ and $%
\left( \ref{c2}\right) $ by $u\left( t\right) $ and $v\left( t\right) ,$
respectively and subtract from each other, then we have, 
\begin{equation}
v\left( t\right) \nabla _{a}^{\mu }\left( p\left( t\right) _{b}\nabla ^{\mu
}u\left( t\right) \right) -u\left( t\right) \nabla _{a}^{\mu }\left( p\left(
t\right) _{b}\nabla ^{\mu }v\left( t\right) \right) =\left[ k\left( t\right)
-m\left( t\right) \right] u\left( t\right) v\left( t\right) .  \tag{15}
\label{c3}
\end{equation}%
Let $t_{1}$ and $t_{2}$ are two succesive zeros of $u.$ If we take the sum
of $\left( \ref{c3}\right) $ from $t_{1}$ to $t_{2}$, then we have%
\begin{equation*}
\sum_{s=t_{1}}^{t_{2}}v\left( s\right) \nabla _{a}^{\mu }\left( p\left(
s\right) _{b}\nabla ^{\mu }u\left( s\right) \right)
-\sum_{s=t_{1}}^{t_{2}}u\left( s\right) \nabla _{a}^{\mu }\left( p\left(
s\right) _{b}\nabla ^{\mu }v\left( s\right) \right) =\sum_{s=t_{1}}^{t_{2}} 
\left[ k\left( s\right) -m\left( s\right) \right] u\left( s\right) v\left(
s\right) ,
\end{equation*}%
and if we apply the integration by parts formula, note that integration by
parts is slightly different in here but it is seen by interchanging order of
summation, then the left hand side of the last equation is zero, so%
\begin{equation*}
\sum_{s=t_{1}}^{t_{2}}\left[ k\left( s\right) -m\left( s\right) \right]
u\left( s\right) v\left( s\right) =0,
\end{equation*}%
this implies that%
\begin{equation*}
k\left( t\right) =m\left( t\right) ,
\end{equation*}%
and hence, we arrive at a contradiction. The theorem is proved.

\bigskip

\noindent \textbf{Theorem 3.2. }(Sturm's 2nd comparison theorem) Let $%
k\left( t\right) <m\left( t\right) $ over the entire interval $\left[ a,b%
\right] ,$ if $u\left( t\right) $ has $n$ zeros in the interval $\left( a,b%
\right] ,$ then $v\left( t\right) $ has not less than $n$ zeros in the same
interval and $k$th zero of $v\left( t\right) $ is less than $k$th zero of $%
u\left( t\right) .$

\bigskip

\textbf{Proof. }Let $t_{1}$ is the zero of $u\left( t\right) $ closest to
the point $a,$ note that suppose $u\left( a\right) \neq 0$. Let's prove that 
$v\left( r\right) $ has at least one zero in the interval $\left(
a,t_{1}\right) $. If we multiply the equation $\left( \ref{c1}\right) $ and $%
\left( \ref{c2}\right) $ by $u\left( t\right) $ and $v\left( t\right) ,$
respectively, subtract from each other and take the sum of the result from $%
a $ to $t_{1}$, then we have by the help of integration by parts%
\begin{equation*}
\sum_{s=a}^{t_{1}}\left[ k\left( s\right) -m\left( s\right) \right] u\left(
s\right) v\left( s\right) =0,
\end{equation*}%
this implies that%
\begin{equation*}
k\left( t\right) =m\left( t\right) ,
\end{equation*}%
and hence, we arrive at a contradiction. The theorem is proved.

\bigskip

Let's consider the DFSL equation secondly in (\textbf{G-L)} sense; $0<\mu <1$%
\begin{eqnarray}
\Delta _{-}^{\mu }\left( \tilde{p}\left( t\right) \Delta _{+}^{\mu }\tilde{u}%
\left( t\right) \right) +\tilde{k}\left( t\right) \tilde{u}\left( t\right)
&=&0,  \TCItag{16}  \label{c4} \\
\Delta _{-}^{\mu }\left( \tilde{p}\left( t\right) \Delta _{+}^{\mu }\tilde{v}%
\left( t\right) \right) +\tilde{m}\left( t\right) \tilde{v}\left( t\right)
&=&0.  \TCItag{17}  \label{c5}
\end{eqnarray}

\bigskip

\noindent \textbf{Theorem 3.3. }(Sturm's 1st comparison theorem) Let $\tilde{%
k}\left( t\right) <\tilde{m}\left( t\right) $ over the entire interval $%
\left[ a,b\right] ,$ then between every two zeros of any nontrivial solution
of the equation $\left( \ref{c4}\right) $ there is at least one zero of
every solution of the equation $\left( \ref{c5}\right) $.

\bigskip

\textbf{Proof. }Proof is similar to the proof of Theorem 3.1.

\bigskip

\noindent \textbf{Theorem 3.4. }(Sturm's 2nd comparison theorem) Let $\tilde{%
k}\left( t\right) <\tilde{m}\left( t\right) $ over the entire interval $%
\left[ a,b\right] ,$ if $\tilde{u}\left( t\right) $ has $n$ zeros in the
interval $\left( a,b\right] ,$ then $\tilde{v}\left( t\right) $ has not less
than $n$ zeros in the same interval and $k$th zero of $\tilde{v}\left(
t\right) $ is less than $k$th zero of $\tilde{u}\left( t\right) .$

\bigskip

\textbf{Proof.}Proof is similar to the proof of Theorem 3.2.

\bigskip

A symbolic graph is given in what follows to clarify the comparison theorem;
the zeros of two linearly independent solutions of the Airy equation $%
y^{\prime \prime }-xy=0$ alternate, as asserted by the Sturm comparison
theorem%
\begin{equation*}
\underset{}{\FRAME{itbpF}{3.3806in}{2.098in}{0in}{}{}{Figure}{\special%
{language "Scientific Word";type "GRAPHIC";maintain-aspect-ratio
TRUE;display "USEDEF";valid_file "T";width 3.3806in;height 2.098in;depth
0in;original-width 2.5469in;original-height 1.5705in;cropleft "0";croptop
"1";cropright "1";cropbottom "0";tempfilename
'P0OW4200.wmf';tempfile-properties "XPR";}}}
\end{equation*}%
\textbf{4. Conclusion}

In this study, we give the Sturm comparison theorems for discrete fractional
Sturm-Liouville (DFSL) equations within Riemann-Liouville and Gr\"{u}%
nwald-Letnikov sense. One of the most important results is Sturm comparison
theorems \cite{levi}. These theorems give information about the properties
of zeros of two equations having different potentials. Besides, these
theorems will be the basis of "oscillation theorems", which has great
importance in spectral theory. This study would be benefit for the theory of
DFSL equations.

\end{document}